\documentclass{proc-l}

\newtheorem{theorem}{Theorem}[section]
\newtheorem{lemma}[theorem]{Lemma}

\theoremstyle{definition}

\theoremstyle{remark}

\numberwithin{equation}{section}

\begin{document}

\title{Fourier Analysis and the closed form for the Zeta Function at even positive integers}


\author{Jibran Iqbal Shah }
\address{Student, Arrowad International School, Riyadh, Saudi Arabia}
\curraddr{}
\email{jibraniqbal2003@gmail.com}
\thanks{}

\keywords{Zeta Function, Number Theory, Analytic Number Theory, Bernoulli Numbers}

\subjclass[2010]{ 11R42 (Primary), 11B68 (Secondary)}

\date{}

\dedicatory{}


\begin{abstract}
Using a summation identity obtained for the Fourier coefficients of $x^{2k}$, we derive a closed form expression for the zeta function at even positive integers, using a technique similar to one in an existing proof by Aladdi and Defant[1], but in a simpler and shorter way.
\end{abstract}

\maketitle


\section*{Introduction}

The Basel problem, more specifically the problem of evaluating the sum
$$\sum_{n=1}^{\infty} \frac{1}{n^{2}}$$
was a very famous problem back in the day. The problem was  popularized by the Bernoullis who lived in Basel, Switzerland, so it became known as the Basel Problem. Leonhard Euler solved it in 1735 and along with solving the Basel Problem, he also found a closed-form evaluation of
$$\zeta(2 k)=\sum_{n=1}^{\infty} \frac{1}{n^{2 k}}$$
for all even integers $2 k \geq 2$. The value of $\zeta(2 k)$ is a rational multiple of $\pi^{2 k}$. Euler showed that
$$\zeta(2 k)=\frac{(-1)^{k+1} B_{2 k}(2 \pi)^{2 k}}{2(2 k) !}$$
Where $B_k$ are the Bernoulli Numbers, note that these are defined as [7]
$$ \frac{t}{e^{t}-1}=\frac{t}{2}\left(\operatorname{coth} \frac{t}{2}-1\right) \quad=\sum_{m=0}^{\infty} \frac{B_{m}^{-} t^{m}}{m !}  $$
$$ \frac{t}{1-e^{-t}}=\frac{t}{2}\left(\operatorname{coth} \frac{t}{2}+1\right) \quad=\sum_{m=0}^{\infty} \frac{B_{m}^{+} t^{m}}{m !} $$
The only difference between $B_m^+$ and $B_m^-$ being that $B_1^+ = \frac{1}{2}$ and $B_1^- = -\frac{1}{2} $. There are many ways to prove Euler's closed form for $\zeta(2k)$, some simpler than others but they appear in [3] and [4].
\\ \\

In this paper, we will, using the $2\pi$ periodic Fourier Series representation of $x^{2k}$ in the interval $[-\pi,\pi]$, derive a recurrence relation for $\zeta(2k)$ and obtain the closed form expression Euler obtained with the help of Bernoulli numbers and a new identity with them. This proof, although uses the same core ideas, takes a different route than that of Alladi and Defant[1] as it avoids the use of Parseval's theorem, and hence ends up being shorter, easier to understand and more elegant. The recurrence for $\zeta(2k)$ in this paper relies on $\zeta(2j)$ for $j < k$, but a paper by Kuo[5] provides a recurrence for $\zeta(2k)$ for $j \leq \frac{n}{2}$. The novelty in this paper is the use of the Fourier Coefficients of $x^{2k}$ rather than $x^k$ and a new Bernoulli Identity which we prove and use to derive our final result. An example of the case for $k = 1$ is given below.

Consider the Fourier Series representation of $f(x) = x^2 $ for $x \in [-\pi, \pi]$
$$ f(x) =  a_0 + \sum_{n=1}^{\infty}a_n\cos(nx) + b_n\sin(nx) $$
$$ a_0 = \frac{1}{2\pi} \int_{-\pi}^{\pi} x^2 \,\textmd{d}x =  \frac{\pi^2}{3} ,\quad
 a_n = \frac{1}{\pi} \int_{-\pi}^{\pi} x^2\cos(nx) \,\textmd{d}x = \frac{4 (-1)^n}{n^2}, \quad
 b_n = 0 $$
$$ f(x) = \frac{\pi^2}{3} + \sum_{n=1}^{\infty} \frac{4 (-1)^n}{n^2}\cos(nx)$$
At $x =\pi$, we obtain 
$$ \pi^2 = \frac{\pi^2}{3} + 4\sum_{n=1}^{\infty} \frac{1}{n^2} \Rightarrow \zeta(2) = \frac{\pi^2}{6} $$

\section{The proof of the generalized sum $\zeta(2k)$}

\begin{lemma} 
$$\sum_{n=1}^{\infty}(-1)^{n} I(n, k)=\frac{k \pi^{2 k+1}}{2 k+1}  \quad \textmd{where}\quad I(n,k) = \int_{0}^{\pi} x^{2k}\cos(nx) \,\textmd{d}x $$
\end{lemma}

\begin{proof}
Consider the Fourier series representation of $f(x) = x^{2k}$ for $k \in \mathbb{N}$ in the interval $[-\pi, \pi]$ and calculate $a_0$, $a_n$ and $b_n$
$$ a_0 = \frac{1}{2\pi} \int_{-\pi}^{\pi} x^{2k} = \frac{\pi^{2k}}{2k+1} , \quad a_n = \frac{1}{\pi} \int_{-\pi}^{\pi}  x^{2k}\cos(nx) \,\textmd{d}x = \frac{2}{\pi} I(n,k) $$
as $x^{2k}\cos(nx)$ is an even function. 
$$ b_n = \frac{1}{\pi}\int_{-\pi}^{\pi} x^{2k}\sin(nx) \,\textmd{d}x = 0 $$
as $x^{2k}\sin(nx)$ is an odd function. Substituting $x = \pi$ yields 
$$ \pi^{2k} = \frac{\pi^{2k}}{2k+1} + \frac{2}{\pi}\sum_{n=1}^{\infty} (-1)^n I(n,k) \Rightarrow \sum_{n=1}^{\infty}(-1)^nI(n,k) = \frac{k\pi^{2k+1}}{2k+1}.$$
\end{proof}

\begin{lemma}
$$I(n, k)=\sum_{i=1}^{k} \frac{(-1)^{n}}{n^{2 i}} \frac{(-1)^{i-1}(2 k) !}{[2 k-(2 i-1)] !} \pi^{2 k-(2 i-1)} \quad \textmd{where} \quad I(n,k) = \int_{0}^{\pi} x^{2k}\cos(nx) \,\textmd{d}x  $$
\end{lemma}

\begin{proof}
$$ I(n,k) = \int_0^{\pi} x^{2k}\cos(nx) \,\textmd{d}x  $$
Carrying out integration by parts repeatedly, one obtains 
$$ I(n,k) =  \left. \frac{(2k)!x^{2k-1} \cos(nx)}{(2k-1)!n^2} - \frac{(2k)!x^{2k-3}\cos(nx)}{(2k-3)!n^4}  + \cdots (-1)^{k-1}\frac{(2k)! x \cos(nx)}{n^{2k}} \right]_0^{\pi} $$
Substituting the limits of integration, we obtain 
$$I(n, k)=\sum_{i=1}^{k} \frac{(-1)^{n}}{n^{2 i}} \frac{(-1)^{i-1}(2 k) !}{[2 k-(2 i-1)] !} \pi^{2 k-(2 i-1)}.$$
\end{proof}

\begin{lemma}
$$\sum_{i=1}^{k} \zeta(2 i)\pi^{2 k-2 i} \frac{(-1)^{i-1}(2 k) !}{[2 k-(2 i-1) !} =\frac{k \pi^{2 k} }{2 k+1}$$
\end{lemma}

\begin{proof}
Substituting the expression obtained for $I(n,k)$ in Lemma 1.2 into Lemma 1.1 (and noting the $(-1)^n$ in Lemma 1.1 canceling out the $(-1)^n$ in Lemma 1.2) we obtain
$$\sum_{n=1}^{\infty} \sum_{i=1}^{k} \frac{(2k)! \pi^{2 k-(2 i-1)} (-1)^{i-1}}{ n^{2i}[2 k-(2 i-1)] ! }=\frac{k \pi^{2 k+1}}{2 k+1}$$
Taking $n^{2i}$ out of the first summation as it is independent of $k$, we obtain
$$ \sum_{i=1}^{k} \zeta(2i) \frac{(-1)^{i-1} (2k)! \pi^{2k-2i+1}}{[2k-(2i-1)]!} = \frac{k\pi^{2k+1}}{2k+1} $$
This is the summation identity mentioned in the abstract. It can be rearranged to provide a recurrence relation for $\zeta(2k)$ in terms of $\zeta(2j)$ where $j<n$. \\ \\ Dividing both sides of the equation by $\pi$ proves the lemma 
$$\sum_{i=1}^{k} \zeta(2 i)\pi^{2 k-2 i} \frac{(-1)^{i-1} (2 k) !}{[2 k-(2 i-1)] !} =\frac{k \pi^{2 k} }{2 k+1}.$$
\end{proof}

Notice that the expression on the right hand side is a single fraction containing a power of $\pi$. This means the terms in the summation all must be multiples of the same power of $\pi$. If that is the case, then each term in the summation must contain $\pi^{2k}$ and hence $\zeta(2i) = C_i \pi^{2i}$ where $C_i$ are coefficients to be determined.

\begin{lemma}
$$\sum_{i=1}^{k} 2^{2 i} B_{2 i}\binom{2k}{2i} \frac{1}{2 k-2 i+1}=\frac{2 k}{2 k+1}$$
\end{lemma}

\begin{proof}

First we prove a similar identity, specifically 
\begin{equation}
    \sum_{i=1}^{m} 2^{i} B_{i}^{+}\binom{m}{i} \frac{1}{m-i+1}=\frac{2 m+1}{m+1}
\end{equation}
Using the property of Bernoulli Polynomials, which are defined as 
$$ B_{n}(x)=\sum_{k=0}^{n}\binom{n}{k} B_{k}^- x^{n-k} $$
More on the definition of the Bernoulli Polynomials and its generalization is present in [6]. We then have 
$$\sum_{i=0}^{m} B_{i}^- \binom{m}{k}\frac{2^{i}}{m-i+1} =2^{m+ 1} \sum_{i=0}^{m} B_{i}^- \binom{m}{i} \int_{0}^{\frac{1}{2}} x^{m-i} \mathrm{~d} x =2^{m+1} \int_{0}^{\frac{1}{2}} B_{m}(x) \mathrm{d} x$$
use the substitution 
$$ y = 1-x \quad \textmd{d}y = -\textmd{d}x $$
$$\int_{0}^{\frac{1}{2}} B_{m}(x) \mathrm{d} x=\int_{\frac{1}{2}}^{1} B_{m}(1-x) \mathrm{d} x=(-1)^{m} \int_{\frac{1}{2}}^{1} B_{m}(x) \mathrm{d} x=\int_{\frac{1}{2}}^{1} B_{m}(x) \mathrm{d} x$$
$(-1)^m = 1$ as $m$ is even and hence using a well known property of Bernoulli Polynomials [2]
$$\int_{0}^{\frac{1}{2}} B_{m}(x) \mathrm{d} x=\frac{1}{2} \int_{0}^{1} B_{m}(x) \mathrm{d} x=0$$
The only difference between $B_1^+$ and $B_1^-$ is that $B_1^{+} = \frac{1}{2}$ and $ B_1^- = -\frac{1}{2} $ Plugging in $i = 0$, subtracting the term with $i=1$ for $B_i^-$ and adding the term for $B_i^+$, and rearranging, we get 
$$ \sum_{i=1}^{m} 2^{i} B_{i}^{+}\binom{m}{i} \frac{1}{m-i+1}=\frac{2 m+1}{m+1} $$

when $m$ is a positive even integer. This can be rewritten as 
\begin{equation}
     \sum_{i=1}^{2k} 2^i B_i^+ \binom{2k+1}{i} = 4k+1 
\end{equation}

If $m$ is an even positive integer, then $m = 2k$ for some $ k \in \mathbb{N}$, this can be rewritten as 
$$ \sum_{i=1}^{2k} 2^i B_i^+ \binom{2k}{i} \frac{1}{2k-i+1} = \frac{2k}{2k+1} + 1 $$
Notice how $\frac{4k+1}{2k+1} = \frac{2k}{2k+1} + 1$. \\ \\ If Eq. (1.1) is true then as $B_{k} = 0$ when $k$ is odd, we have 
$$  \sum_{i=1}^{2k} 2^i B_i^+ \binom{2k}{i} \frac{1}{2k-i+1} -1  =  \sum_{i=1}^{k} 2^{2 i} B_{2 i} \binom{2k}{2i} \cdot \frac{1}{2 k-2 i+1}  $$
The $-1$ is present because at $i = 1$, the expression evaluates to 1 ($B_1 =\frac{1}{2}$). This completes the proof of the identity
$$ \sum_{i=1}^{k} 2^{2 i} B_{2 i} \binom{2k}{2i} \frac{1}{2 k-2 i+1}  = \frac{2k}{2k+1}. $$

\end{proof}

\begin{lemma}
$$ \textmd{If} \quad \sum_{i=1}^{k} C_{i} \frac{(-1)^{i-1}(2 k) !}{(2 k-2 i+1) !}=\frac{k}{2 k+1} \quad \textmd{then} \quad C_i = \frac{(-1)^{i-1} 2^{2 i} B_{2 i}}{2(2 i) !}$$
\end{lemma}

\begin{proof}
Begin with Lemma 1.4 which states
$$ \sum_{i=1}^{k} 2^{2 i} B_{2 i}\binom{2k}{2i}\frac{1}{2 k-2 i+1}=\frac{2 k}{2 k+1} $$ 
Substituting in the identity, 
$$ \binom{2k}{2i} \frac{1}{2k-2i+1} = \frac{(2k)!}{(2i)!(2k-2i+1)!} $$
We obtain 
 $$ \sum_{i=1}^{k} 2^{2i-1}B_{2i} \frac{(2k)!}{(2i)!(2k-2i+1)!} = \frac{k}{2k+1} $$
 As both expressions sum up to $\frac{k}{2k+1}$, we can equate them to give 
 $$ \sum_{i=1}^{k} 2^{2i-1}B_{2i} \frac{(2k)!}{(2i)!(2k-2i+1)!} = \sum_{i=1}^{k} C_{i} \frac{(-1)^{i-1}(2 k) !}{(2 k-2 i+1) !}  $$
Comparing both sides, we obtain
 $$ C_i = \frac{(-1)^{i+1} 2^{2 i} B_{2 i}}{2(2 i) !}. $$
\end{proof}

 \begin{theorem}
 $$ \zeta(2i) = 
(-1)^{i+1} \frac{(2 \pi)^{2 i} B_{2 i}}{2(2 i) !}
 $$
For an even positive integer $i$
 \end{theorem}

  \begin{proof}
  Substituting $\zeta(2i) = C_i\pi^{2i}$ into Lemma 1.3 we get 
  $$  \sum_{i=1}^{k} C_i\pi^{2 k} \frac{(-1)^{i-1}(2 k) !}{[2 k-(2 i-1)] !} =\frac{k \pi^{2 k} }{2 k+1} $$
  Dividing both sides by $\pi^{2k}$ we obtain : 
  $$  \sum_{i=1}^{k} C_i \frac{(-1)^{i-1}(2 k) !}{(2 k-2 i+1) !}=\frac{k }{2 k+1}  $$
 As $\zeta(2i) = C_i \pi^{2i}$, according to Lemma 1.5,  we obtain 
  $$ \zeta(2i) = \frac{(-1)^{i+1} (2\pi)^{2i} B_{2 i}}{2(2 i) !}. $$
 \end{proof}

 \section{Remarks}
In this paper we obtained and used the new identity
$$ \sum_{i=1}^{k} 2^{2 i} B_{2 i}\binom{2k}{2i} \frac{1}{2 k-2 i+1}=\frac{2 k}{2 k+1}  $$
 In the proof of Lemma 1.4, we have Eq. (1.2) which stated that 
$$ \sum_{i=1}^{2k} 2^i B_i^+ \binom{2k+1}{i} = 4k+1 $$
Readers can attempt to generalize the sum by considering finding a closed form for
$$ \sum_{i=1}^{2kx} x^i B_i \binom{2kx+1}{i}. $$

\section{Acknowledgements}
I would like to thank Colin Defant, Sai Teja Somu, Juan Luis Varona and Fabio M. S. Lima for spotting typographical errors and providing me with their insightful comments on the paper

\bibliographystyle{spr-chicago}      
\bibliography{example}   
\nocite{*}


\bibliographystyle{amsplain}

\end{document}